\documentclass[11pt, a4paper, leqno]{article}
\usepackage{amssymb}
\usepackage{amscd}
\usepackage{amsmath}
\usepackage{amsthm}
\usepackage{mathrsfs}
\usepackage{graphics}
\usepackage[dvipdfmx]{graphicx}
\usepackage{color}

\usepackage[all]{xy}
\usepackage{titlesec}

\def\overset#1#2{{\mathrel{\mathop {{#2}_{}}\limits^{#1}}}}
\def\underset#1#2{{\mathrel{\mathop {{}_{} {#2}}\limits_{{#1}_{}}}}}
\def\upplim_#1{\underset{#1}{\overline\lim}\;}
\def\lowlim_#1{\underset{#1}{\underline\lim}\;}

\makeatletter
\def\@makefnmark{\hbox{\@textsuperscript{\normalfont
\@thefnmark)}}}
\renewenvironment{enumerate}
  {\ifnum \@enumdepth >3\relax\@toodeep\else
   \advance\@enumdepth\@ne
   \edef\@enumctr{enum\romannumeral\the\@enumdepth}%
   \list{\csname label\@enumctr\endcsname}{%
         \ifnum \@listdepth=\@ne \topsep.1\normalbaselineskip
           \else\topsep\z@\fi
         \parskip\z@ \itemsep\z@ \parsep\z@
         \labelwidth1pc \labelsep0.5pc
         \ifnum \@enumdepth=\@ne \leftmargin1pc\relax
           \else\leftmargin\leftskip\fi
         \advance\leftmargin +1.5pc
         \usecounter{\@enumctr}%
         \def\makelabel##1{\hss\llap{##1}}}%
   \fi}{\endlist}
\makeatother
\newtheorem{cor}[equation]{Corollary}
\newtheorem{defn}[equation]{\indent{\it Definition}\rm }

\newtheorem{lem}[equation]{Lemma}
\newtheorem{prop}[equation]{Proposition}
\newtheorem{rmk}[equation]{\indent \rm {\it Remark}}

\newtheorem{thm}[equation]{Theorem}

\newtheorem{probl}[equation]{\indent Problem}
\newtheorem{quest}[equation]{\indent Question}

\newcommand{\sA}{\mathscr{A}}
\newcommand{\C}{{\mathbf{C}}}

\newcommand{\sD}{\mathscr{D}}

\newcommand{\del}{{\partial}}
\newcommand{\delbar}{\bar{\partial}}

\newcommand{\cE}{\mathcal{E}}

\newcommand{\bG}{\mathbf{G}}
\newcommand{\sG}{\mathscr{G}}
\newcommand{\bH}{\mathbf{H}}
\newcommand{\iso}{\cong}

\newcommand{\im}{{\Im\,}}
\newcommand{\Ker}{{\mathrm{Ker}\,}}

\newcommand{\sL}{\mathscr{L}}
\newcommand{\Lie}{{\mathop{\mathrm{Lie}}}}
\newcommand{\lto}{\longrightarrow}

\newcommand{\N}{\mathbf{N}}

\renewcommand{\O}{{\mathcal{O}}}

\renewcommand{\P}{{\mathbf{P}}}
\newcommand{\rP}{{\mathrm{P}}}

\newcommand{\pNc}{{\mathbf{P}^N(\mathbf{C})}}
\newcommand{\ponec}{{\mathbf{P}^1(\mathbf{C})}}
\newcommand{\ptwoc}{{\mathbf{P}^2(\mathbf{C})}}

\newcommand{\R}{{\mathbf{R}}}

\newcommand{\re}{{\Re\,}}
\newcommand{\sS}{\mathscr{S}}

\newcommand{\supp}{\mathrm{Supp}\,}
\newcommand{\tensor}{\otimes}

\newcommand{\sX}{\mathscr{X}}

\newcommand{\Z}{\mathbf{Z}}

\newenvironment{pf}{\par{\it Proof.\,}}{\qed\vskip+2pt}
\numberwithin{equation}{section}

\title{Analytic and rational sections of relative\\ semi-abelian varieties}
\author{P. Corvaja
,~ J. Noguchi\footnote{\,Supported by
Grant-in-Aid for Scientific Research (C) 19K03511.}
~ and U. Zannier
} 
\begin{document}
\setlength{\baselineskip}{15pt}
\parskip+2.5pt
\maketitle
\thispagestyle{empty}
\begin{abstract}
The hyperbolicity statements for subvarieties and complement of
hypersurfaces in abelian varieties admit arithmetic analogues, due to
Faltings, Ann.\ Math.\ {\bf 133}  (1991)
 (and for the semi-abelian case, Vojta,
 Invent.\ Math.\ {\bf 126} (1996);
 Amer.\ J.\ Math.\ {\bf 121} (1999)).
 In Atti Accad.\ Naz.\ Lincei Rend.\ Lincei Mat.\ Appl.\ {\bf 29} (2018)
 by the second author, an analogy between the analytic and arithmetic
 theories  was shown to hold also at proof level, namely in a 
 proof of Raynaud's theorem (Manin-Mumford Conjecture).
The first aim of this paper is to extend to the relative setting the
above mentioned hyperbolicity results. We shall be concerned with
analytic sections of a relative (semi-)abelian scheme $\sA \to B$ over
 an affine algebraic curve $B$.
 These sections form a group;
while the group of the rational sections (the Mordell-Weil group)
 has been widely studied,  little investigation has been pursued so far
 on the group of the analytic sections.
We take the opportunity of developing some basic structure of this
apparently new theory, defining a notion of height or order
functions for the analytic sections, by means of Nevanlinna theory.
\end{abstract}



\section{Introduction}\label{intr}
 Let $A$ be an { abelian} variety and let $f:\C\to A$
be an entire curve on it. Then the Zariski-closure of its image is a
translate of an abelian subvariety of $A$ (Bloch-Ochiai's Theorem,
c.f., e.g., \cite{nw14} \S4.8), and the same  holds
 for semi-abelian varieties (cf.\ ibid.).

In another direction, it is known that
if $D\subset A$ is an ample divisor of an abelian variety, there is
no non-constant entire curve on $A\setminus D$.
 Analogous results hold for semi-abelian varieties (cf., e.g., \cite{nw14} Chap.~6).
See a relevant Question \ref{top} about the
 {\it topological} closure of an entire curve.

These hyperbolicity statements for subvarieties and complements of
hypersurfaces in abelian varieties admit arithmetic analogues, due to
Faltings \cite{falt} (and Vojta \cite{vo} for the semi-abelian case): the rational
points on a subvariety of an abelian variety are contained in a finite
union of translates of  abelian subvarieties and the integral points on
the complement of an ample divisor in an abelian variety  are finite in
number. 
 In \cite{n18}, the second author observed a direct
relation between those analytic results and  Diophantine properties
in a proof of M. Raynaud's Theorem (Manin-Mumford's Conjecture),
 going beyond formal analogies holding at the level of statements.   

The first aim of this paper is to extend to the relative setting the
above (nowadays classical) analytic results; in our present situation,
the single abelian variety $A$ will be replaced by an algebraic family
$\sA\to B$ of abelian varieties over an algebraic base $B$, and we shall
be concerned with 
 {\it possibly transcendental}
holomorphic sections $B\to \sA$. These sections form a group,
which is the complex-analytic analogue of the group of rational
sections; while this last group, called the Mordell-Weil group of the
abelian scheme $\sA\to B$, has been widely studied, little investigation
has been pursued so far on the group of analytic sections. In this work
we take the opportunity of developing some basic structure of this
apparently new theory, defining e.g. a notion of height for the
transcendental sections, by means of Nevanlinna theory.
\smallskip

We start by noticing the following:  Let $\pi: \sA \to B$ be a holomorphic family
of principally polarized abelian varieties over a base space
$B$ which is algebraic. Then the family $\pi: \sA \to B$ is
algebraic. This is a  result by   Kobayashi--Ochiai \cite{KO71},
 in the spirit of a `Big Picard  Theorem'.
 Hence, the initial datum will consist in an algebraic
family, but we shall consider (possibly) transcendental  sections.

More precisely, our main concern will be addressed to the following
three problems below, motivated by the known results in the constant case:
Let $\pi:\sA \to B$ be an algebraic family of semi-abelian varieties.
In this paper we always assume the existence of a section $B\to \sA$,
namely the  zero-section.
 Then:
 
\begin{probl}\label{st1}
Let $\sD$ be a relatively big divisor on $\sA$ over $B$.
Then, every holomorphic section $\sigma: B \to \sA\setminus \sD$
(omitting $\sD$) is rational.
\end{probl}

\begin{probl}\label{st2}
Let $\sigma:B \to \sA$ be a transcendental holomorphic section,
and let $X(\sigma)$ be the Zariski-closure of $\sigma(B)$ in $\sA$.
Then, $X(\sigma)$ contains a translate of a relative non-trivial subgroup
of $\sA$ over $B$.
\end{probl}

The conclusion that $X(\sigma)$ is itself a translate of a relative
subgroup cannot always hold (see \S\ref{transc} (a)).
However, we may conjecture that
this is the case whenever the section satisfies a stronger condition, to
be precisely formulated, named strict transcendency (see Definition \ref{strtr}):

\begin{probl}\label{st3}
Let $\sigma:B \to \sA$ be a strictly transcendental holomorphic section,
and let $X(\sigma)$ be the Zariski-closure of $\sigma(B)$ in $\sA$.
Then, $X(\sigma)$ is a translate of a relative subgroup
of $\sA$ over $B$.
\end{probl}

We shall start from the concrete example of the Legendre elliptic scheme:
\begin{equation*}
ZY^2=X(X-Z)(X-\lambda Z) \subset B\times \P^2 ,
\end{equation*}
where $\lambda$ varies on the base
$B=\P^1\setminus\{0,1,\infty\}=\C\setminus\{0,1\}$.
 For each $\lambda\in B$ the above curve $\cE_\lambda$, together with its distinguished
point $(0:1:0)$, is an elliptic curve.  
Removing this point from each fiber, so removing a relatively ample
divisor, we obtain the affine variety of equation
 \begin{equation*}
 y^2=x(x-1)(x-\lambda),
 \end{equation*} 
still fibered over the curve $B$. 
As an example of a result in the direction of Problem \ref{st1}, we
shall prove, by means of Yamanoi's Second Main Theorem \cite{y06}, that
any holomorphic section $B\ni \lambda\mapsto (x(\lambda),y(\lambda))$ 
of this fibration is reduced to one of the three
 $2$-torsion points with $y(\lambda)= 0$
(see Theorem \ref{legendre}); in the course of the proof
we show a rationality criterion of holomorphic sections for a
base-extended family of the Legendre elliptic scheme (see Theorem \ref{rat})

We will also see a similar property for hyperelliptic schemes of higher
genera by making use of an extension theorem of big Picard type
due to Noguchi \cite{n85}.

A general result of that type for sections of a family of semi-abelian varieties
will be proved in Theorem \ref{bp}, under the stronger hypothesis that
the family admits no bad reduction. Note that this hypothesis excludes
the non-isotrivial elliptic schemes.

To prove our main results, we will generalize the Nevanlinna theory of
holomorphic curves
in semi-abelian varieties (\cite{nw14}) to a relative setting,
and prove a Big Picard Theorem for a local smooth family;
in this context, we have an interesting problem when the family is singular
over the special point.

\section{Sections of Legendre scheme and related examples.}
Let $(X:Y:Z)\in \ptwoc$ denote a
homogeneous coordinate system of $\ptwoc$.
We consider the simplest but fundamental  Legendre scheme:
\begin{align}\label{leg}
B &= \C \setminus \{0,1\},\\ \notag
\cE &= \{(\lambda, (X:Y:Z))\in B \times \ptwoc :
 Y^2 Z=X(X-Z)(X-\lambda Z) \}, \\
\pi &: \cE \ni (\lambda, (X:Y:Z)) ~ {\mapsto}  ~
 \lambda \in  B, \notag\\
\cE_\lambda &= \pi^{-1}\{\lambda\}. \notag
\end{align}
We set
\[
 P_1=(0: 0: 1), ~ P_2=(1:0:1), ~ Q=(0:1:0).
\]
We can view $\cE$ as a hypersurface in $B \times \ptwoc$;
 the natural projection $\pi:\cE \to B$ together with the (zero) section
$B \times \{Q\}$ gives it the structure of an elliptic scheme over $B$.
We call it the Legendre scheme. It admits three sections of order 2,
 namely
\[
 B \times \{P_1\} , \quad  B \times \{P_2\},
\]
which are abbreviated as $P_1, P_2$,
and their sum
\[
P_3:= (P_1 + P_2) :   B\ni \lambda  \mapsto (\lambda, (\lambda  : 0 : 1))
\in \cE \setminus Q,
\]
where $Q$  denotes also the section $B \times \{Q\}$.

\subsection{Examples}
Before presenting a theorem in the particular case of the Legendre scheme,
we discuss examples of elliptic schemes and their sections.
We keep the notation given above.

\subsubsection{Rational sections omitting a relatively ample divisor.}
 We show in this sub-section that rational sections of abelian schemes
 can indeed omit relatively ample divisors in a non-trivial way.
We first construct a non-isotrivial algebraic family of abelian varieties
over an algebraic variety. 
Set
\begin{align*}
E_\lambda &= \cE_\lambda \setminus\{Q\}=\{y^2=x(x-1)(x-\lambda)\} \subset \C^2,
 \quad \lambda \in B,\\
\tilde B &=\bigcup_{\lambda \in B} \{(\lambda, (x,y)): (x,y) \in E_\lambda\}
\subset B \times \C^2, \\
\sA &= \bigcup_{(\lambda, (x,y)) \in \tilde B} \{(\lambda, (x,y),
(u,v)): (u,v) \in \cE_\lambda\} \subset \tilde B \times \P^2 .
\end{align*}
Now, let $\varpi :\sA \to \tilde B$ be the natural projection.  
In this example, 
there are three sections coming from $P_j, 1 \leq j \leq 3$; i.e.,
\begin{align*}
(\lambda, (x,y)) \in \tilde B &\lto  (\lambda, (x,y), (0,0)) \in \sA, \\
(\lambda, (x,y)) \in \tilde B &\lto  (\lambda, (x,y), (1,0)) \in \sA, \\
(\lambda, (x,y)) \in \tilde B &\lto  (\lambda, (x,y), (\lambda,0)) \in
 \sA.
\end{align*}
These omit a relatively ample divisor $\sD$ defined by
\[
 \sD =\tilde B \times \{Q\} \subset \sA.
\]
Other than these, we have
\begin{align*}
& \tau : t=(\lambda, (x,y)) \in \tilde B \lto (\lambda, (x,y), (x,y)) \in \sA,
\\
&\tau(\tilde B) \cap \sD =\emptyset.
\end{align*}

Note that $\tau$ is ``non-constant'' (the definition is subtle). Note
also that by cutting $\tilde B$ we can produce examples where the base
is an affine curve, all its points
at infinity are points of bad reduction for the elliptic scheme
and some rational section omits the divisor at infinity. These examples
include  the so-called Masser's sections, e.g. the section
$\lambda\mapsto (\lambda;,2,\sqrt{2(2-\lambda)})$ on a (ramified) base
change of the Legendre scheme.

\subsubsection{Transcendental sections.}
 Given an elliptic curve $\cE_\lambda$ in the Legendre family, the
 elliptic exponential is well defined as a map
 $\mathrm{Lie}(\cE_\lambda) \iso \C\to \cE_\lambda$. As we shall
 explain in section \ref{ell-sch}, we can identify globally the line
 bundle of the Lie algebras $ \mathrm{Lie}(\cE_\lambda)$, for $\lambda
 \in \tilde{B}$, with the trivial bundle $\tilde{B}\times \C$. Hence we
 shall view the exponential map as a map 
\[
 \exp_\lambda : \C \lto  \cE_\lambda,
\]
 and set
\[
 \sigma: t=(\lambda, (x,y)) \in \tilde B \lto
(\lambda, (x,y), \exp_\lambda (\varphi(x))) \in \sA_t
:=\pi^{-1}\{t\} \subset \sA,
\]
where $\varphi(x)$ is any non-constant polynomial
(or even, entire function) in $x\, (\in \C)$.
Then, $\sigma$ is a {\em transcendental holomorphic section}
of $\varpi: \sA \to \tilde B$.
In this example, we have that $\sigma(\tilde B) \cap \sD \not=\emptyset$.

\subsubsection{Local sections.}
One may easily obtain at least locally a holomorphic
non-rational section in $\cE\setminus Q \overset{\pi}{\lto} B$
 about the boundary points
of $B$, to say, $\lambda=0$. For example, let $\phi(\lambda)$
be a holomorphic function in a neighborhood of $0$ such that
$\phi(0)\not=0$. With
\[
 x=\frac{X}{Z}, \quad  y=\frac{Y}{Z}  \qquad \hbox{in } \cE \setminus Q,
\]
 we then set
\[
 x(\lambda)= \frac{\phi(\lambda)}{\lambda^2}.
\]
We have
\begin{align*}
 y(\lambda)^2 &=\frac{\phi(\lambda)}{\lambda^2}
\left(\frac{\phi(\lambda)}{\lambda^2} -1 \right)
\left(\frac{\phi(\lambda)}{\lambda^2} - \lambda \right)\\
&= \frac{\phi(\lambda)(\phi(\lambda)-\lambda^2)
(\phi(\lambda) - \lambda^3)}{\lambda ^6}.
\end{align*}
Therefore, taking $\delta>0$ small enough, we have
a one-valued branch
\[
y(\lambda) =
 \frac{
\sqrt{\phi(\lambda)(\phi(\lambda)-\lambda^2)
(\phi(\lambda) - \lambda^3)
}
}{\lambda^3}, \qquad 
|\lambda|< \delta.
\]
Even if $\phi(\lambda)$ is a polynomial, $y(\lambda)$ is not
 rational, unless it vanishes identically.

\subsection{Legendre scheme.}\label{legsch}
Globally, we are going to prove:
\begin{thm}\label{legendre}
Let $\cE \overset{\pi}{\lto} B$ be as above in \eqref{leg}.
 Then there is no holomorphic section
$B \to \cE \setminus Q$ other than $P_j ~ (j=1,2,3)$.
\end{thm}

We first prove the rationality of the sections,
in general even after finite base changes (extensions);
 this is the crucial part of the proof, which makes use of 
a deep theorem of Yamanoi, 
and the result may have  an interest of its own.

Let $\phi: \tilde B \to B$ be a finite base change
(i.e., a finite proper rational holomorphic map)
 and let
\begin{equation}\label{extleg}
\tilde\pi:\tilde \cE=\tilde B \times_B \cE \to \tilde B
\end{equation}
be the lift of $\cE/B$. Then $\tilde\cE/\tilde B$ carries the natural
structure of a group scheme induced from $\cE/B$ with the
zero section $\tilde Q$ induced by $Q$.
In general, $\tilde\cE \to \tilde B$ may carry a non-torsion rational section
and hence infinitely many rational sections.

\begin{thm}\label{rat}
Let $\tilde\pi: \tilde\cE \to \tilde B$ be as above in \eqref{extleg} and
let $\gamma: \tilde B \to \tilde\cE$ be a rational section of
 $\tilde\cE \to \tilde B$.
Then, a holomorphic section $\sigma: \tilde B \to \tilde\cE$
is rational if and only if the intersection
 $\sigma(\tilde B) \cap \gamma(\tilde B)$ is finite.
\end{thm}
\begin{pf}
It suffices to prove the ``if'' part.
Replacing $\sigma$ by $\sigma - \gamma$, we may assume that
 $\gamma=\tilde Q$.

We set $f(\lambda)=\lambda(\lambda-1) ~(\lambda \in B)$.
By the embedding 
$$
\lambda \in B \to (\lambda,  1/f(\lambda))
 \in \{(\lambda, \mu): f(\lambda)\mu=1\} \subset  \C^2,
$$
we identify $B$ with the image,
 which is a closed affine algebraic curve in $\C^2$.
We consider $B$ as a ramified cover of $\C$ via
\[
 \pi_B :  B\ni \lambda  \mapsto z=\lambda+ \frac1{f(\lambda)}  \in \C.
\]
We set $\pi_{\tilde B}=\pi_B \circ \phi: \tilde B \to \C$.
Then, $|z|=|\pi_{\tilde B}(\zeta)|$ is an exhaustion function on
 $\tilde B$
by which we define the order function $T_\dagger(r; \star )$ of a meromorphic
functions on $B$, and the counting function $N_k(r, \bullet)$ of a
 divisor truncated at level $k$
 etc.\ (cf.\ \cite{n76}, \cite{nw14} \S3.3.3).

We  compactify
 $\tilde B\hookrightarrow \bar{\tilde B} ~(\subset \pNc)$
with some $\pNc$ and may assume that the polar divisor
$(\phi)_\infty$ of $\phi$ belongs to the
linear system $|\O_{\bar{\tilde B}}(1)|$.
 We set
$X= \bar{\tilde B} \times \ponec$,
 which is provided with the first (resp.\ second) projection
$p: X \to \bar{\tilde B}$ (resp.\  $q: X \to  \ponec$); 
 now every section $\sigma$ determines a  holomorphic map
\[
 g:  \tilde B \ni \zeta  \longmapsto (\zeta, x(\zeta)) \in X.
\]
Using the affine coordinate $w$ of $\C \subset \ponec$,
we regard $w=x(\zeta)$ as a meromorphic function on $\tilde B$
with poles at those $\zeta \in \tilde B$ such that
 $\sigma(\zeta)=\tilde Q(\zeta)$.
Then the section $\sigma$ also provides
 a meromorphic function $y(\zeta)$ on $\tilde B$ satisfying
\begin{equation}\label{legeqext}
 (y(\zeta))^2=x(\zeta)(x(\zeta) -1)(x(\zeta) - \phi(\zeta)).
\end{equation}
 We define the following effective divisor on $X$:
\[
 D=\{w=0\}+\{w=1\}+\{w=\infty\}+\{w-\phi(\zeta)=0\}.
\]
Note that $\pi_{\tilde B}(\zeta)$ 
 and $p \circ g(\zeta)=\zeta$  are rational functions on $\tilde B$.
 So, the order functions satisfy
\[
 T_g(r; L)=T_x(r)+O(\log r),
\]
where $L=p^*\O_{\bar{\tilde B}}(1) \tensor q^*\O_\ponec(1)$
 and $T_x(r)=T_x(r ; \O_\ponec (1))$.
Then we have by Yamanoi \cite{y06} Theorem 1.2 with similar notation
 that
\begin{align*}
 T_g\left(r ; K_{X/\bar{\tilde B}}(D)\right) &\leq {N_1}(r, \{x(\zeta)=0\})+
{N_1}(r, \{x(\zeta)=1\})  + {N_1}(r, \{x(\zeta)=\infty\}) \\ &\quad +
{N_1}(r, \{x(\zeta)-\phi(\zeta)=0\})+ \epsilon T_g(r;  L)+ O(\log r) ||,
\end{align*}
where
$\epsilon$ is an arbitrary small positive number and the symbol $||$
is used for the standard sense in Nevanlinna theory, while  the implicitly
 mentioned exceptional set depends on $\epsilon$.

By the First Main Theorem
\begin{align*}
N(r, \{x(\zeta)=0\}) &\leq T_x(r)+O(1),\\
N(r, \{x(\zeta)=1\}) &\leq T_x(r)+O(1).
\end{align*}
Similarly, since $\{w-\phi(\zeta)=0\}$ is an element of $|L|$,
\[
 N(r, \{x-\lambda=0\}) \leq T_g(r ;  L)+O(1)=T_x(r)+O(\log r),
\]
and since $K_{X/\bar{\tilde B}}=q^*\O_\ponec(-2)$,
\[
 T_g(r ; K_{X/\bar{\tilde B}}(D)))=2 T_x(r)+O(\log r).
\]
By the assumption, $\{\zeta \in \tilde B: \sigma(\zeta)\in \tilde Q\}=
\{\zeta \in \tilde B: g(\zeta)=\infty\}$ is a finite set, so that
\[
 N(r, \{x(\zeta)=\infty\})=O(\log r).
\]
Since all zeros of
 $x(\zeta), x(\zeta)-1, x(\zeta)-\phi(\zeta)$ have order $\geq 2$,
 we have
\begin{align*}
(2-\epsilon)T_x(r)  &\leq  \frac{1}{2} \big( N(r, \{x(\zeta)=0\})
 +N(r, \{x(\zeta)=1\}) + N(r, \{x(\zeta)-\phi(\zeta)=0 ) \big)  
   +O(\log r)||\\
&\leq \frac{3}{2}T_x(r)+O(\log r)||.
\end{align*}
Therefore, by taking $\epsilon < 1/2$
\begin{equation}\label{olog}
T_x(r)=O(\log r)|| ;
\end{equation}
this implies the rationality of
$x(\zeta)$ and hence that of $y(\zeta)$.
\end{pf}

The next lemma finishes the proof of Theorem \ref{legendre}.
\begin{lem}\label{torsion-on-legendre}
The Legendre scheme admits no rational sections other than
$P_j\, (j=1,2,3)$ and $Q$; i.e. the Mordell-Weil group of $\cE \to B$
consists of the $2$-torsion group.
\end{lem}
\begin{pf}
The family of cubic curves $\cE_\lambda=\pi^{-1} (\lambda),~
\lambda \in B \subset \ponec$, forms a pencil having
 as base points the three
points $P_1,  P_2, Q \in \ptwoc$.
Let $\sigma: \lambda \in B \to (\lambda, (X(\lambda) : Y(\lambda :
 Z(\lambda)) \in \cE_\lambda \subset \cE$ be a rational section.
Let $ C \subset \ptwoc$ denote the closure of the projection
of its image $\sigma(B)$ in $\ptwoc$.
Then,  $C$ intersects each curve $E_\lambda$ in three fixed points
$P_1,  P_2, Q$ and possibly a fourth moving point
$\sigma(\lambda) = (X(\lambda) : Y(\lambda) : Z(\lambda))$.
 We shall prove that this point is of $2$-torsion.
 Suppose that $C$ does not reduce to a point (i.e. 
$\sigma(\lambda)$ is not identically equal to $Q$, nor $P_1$ nor 
$P_2$) and let $F(X, Y, Z) = 0$ be an equation for $C$,
 where $F$ is a homogeneous form of degree $d > 0$. Set
\[
 f := \frac{F}{Z^d} \in  \C(\ptwoc).
\]
Consider any value of $\lambda \in  B$ and view $f$
 as a rational function on $\cE_\lambda$.
 The support of the divisor $(f)$ of $f$ is then contained in
$\{Q,  P_1,  P_2, \sigma(\lambda)\}$.
 Identifying $\cE_\lambda$ with its Jacobian
 $\mathrm{Pic}^0(\cE_\lambda)$,
 we obtain that the sum of the elements in $(f)$,
 in the sense of the group law on $\cE_\lambda$, must vanish.
 It follows that $\sigma(\lambda)$ belongs to the
group generated by $P_1,  P_2$, i.e. to the 2-torsion group.
\end{pf}
\begin{rmk}\rm
\begin{enumerate}
\item
We provide here an alternative proof of the vanishing of the Mordell-Weil rank.
 Consider the surface $S$ obtained by blowing-up $\ptwoc$
 over the base locus of the
 pencil of cubics $\cE_\lambda,~ \lambda \in B$,
 i.e. over $P_1,  P_2, Q$ in the above notation.
The surface $S$ is endowed with a natural projection
$\varpi: S  \to \ponec$, which is a well-defined morphism,
whose fibers are the curves $\cE_\lambda$ for $\lambda \in \ponec$.
 By taking for the zero-section the natural map inverting the
$\varpi$ on the exceptional divisor over $Q$
 one obtains a structure of elliptic surface on $S$.
Recall the Shioda-Tate formula (cf.\ Shioda \cite{sh}) for the rank $r$
 of the Mordell-Weil group:
\[
 r = \rho - 2 - \sum_{\lambda \in \ponec} (n_\lambda -1),
\]
where $\rho$ is the Picard number of $S$ and for each point
$\lambda \in \ponec$, $n_\lambda$ is
 the number of components of the fiber $\varpi ^{-1}\lambda$.
 In our case $\rho = 4$ and
 the only reducible fiber is the fiber of $\lambda=\infty$,
 which has three components.
It follows that $r = 0$, i.e. the Mordell-Weil group is torsion.
\item
In the quoted paper,   Shioda proves more: the Mordell-Weil group of
 an elliptic scheme obtained from the Legendre scheme by   any
 {\em unramified} base change is torsion. 
 The interested reader is
 addressed to \cite{cz} for a history, motivations and generalizations
 of this result, as well as different approaches to its proof. 
\end{enumerate}
\end{rmk}

Related to the above problem, N. Katz  asked in a conversation with 
U. Zannier for the case of the hyperelliptic scheme of genus
$ g>1$ of $\P^2_B$  defined by
\begin{equation}\label{hypel}
   y^2 = h(x) (x-\lambda)
\end{equation}
in terms of an affine coordinate system $(x,y) \in \C^2 \subset \ptwoc$,
where $h$ is a given polynomial with complex coefficients
 (i.e., independent of $\lambda$)
 of  degree $ 2g>2$  with simple roots,
and $B=\C\setminus\{h=0\}$. Note that for each $\lambda\in B$ the
equation defines a smooth affine curve with a single point at infinity.
 This family is relevant to Katz' work on monodromy.

In this context we can prove:
\begin{thm} Let $\sX \to B$ be the hyperelliptic scheme defined
by completing  the curves of equation \eqref{hypel} above. Let
$$
\sigma: \lambda \in B \to (\lambda, (x(\lambda), y(\lambda))) \in \sX
$$
be an arbitrary holomorphic section.
 Then, $\sigma$ is a rational
section such that either $y(\lambda)\equiv 0$ or $y(\lambda)\equiv \infty$.
\end{thm}
\begin{pf}
This is a case to which a big Picard theorem (a holomorphic extension
theorem) obtained by Noguchi \cite{n85} is applicable, since $\sX$
is a family of compact curves of genus $\geq 2$; hence, Theorem (5.2)
and Lemma (2.1)   of \cite{n85} imply that  
 $x(\lambda)$ and $y(\lambda)$ are rational functions.

To conclude the proof, we need  to prove  
 an analogue of Lemma \ref{torsion-on-legendre}, namely 
 that:

{\it Claim.
 The only rational points of $\sX$ over $\C(\lambda)$ are those with $y=0$.}

The argument  below is   a variation on the elementary proof of
   ``{\it abc}'' over function fields. 
We may suppose that $h$ is monic and let
 $\xi_1,\ldots ,\xi_{d}\in\C$ be its distinct roots. 

Let $(u(\lambda),v(\lambda))$ be a solution of \eqref{hypel} in rational
 functions,
where $v\neq 0$. This last condition implies that $u, v$ are in fact
 both non-constant. 

Any pole on $\C$ of $u$ or $v$ must appear with even order $2e$  in $u$
 and order  $(d+1)e$ in $v$  for some integer $e$,
and hence we may write 
$$
u(\lambda)={a(\lambda)\over q^2(\lambda)},\qquad
 v(\lambda)={b(\lambda)\over q^{d+1}(\lambda)}
$$
for complex polynomials $a,b$ and $q\neq 0$,
 which are pairwise coprime. This yields the equation
\begin{equation}\label{E.2}
b^2(\lambda)=(a(\lambda)-\lambda q^2(\lambda))
\prod_{i=1}^d(a(\lambda)-\xi_iq^2(\lambda)).
\end{equation}

Since $a,q$ are coprime, the $d$  factors in the product on the right
 are non-zero and pairwise coprime, whereas 
 the $\gcd((a(\lambda)-\lambda q^2(\lambda)),
 a(\lambda)-\xi_i q^2(\lambda))$ divides $\lambda-\xi_i$.
 Hence, since the whole product
 is a (non-zero) square, we must have 
\begin{equation}\label{E.3}
a(\lambda)-\xi_iq^2(\lambda)=(\lambda - \xi_i)^{\mu_i}
c_i^2(\lambda),\qquad  i=1,\ldots ,d,
\end{equation}
for $\mu_i\in\{0,1\}$ and suitable non-zero polynomials $c_i(\lambda)$.
Note that the $c_i(\lambda)$ are pairwise coprime. 

Let $m=\deg u=\max(\deg a,2\deg q)>0$. Note that each factor in the
 product on the right  of \eqref{E.2} (namely, each side of equation
 \eqref{E.3}) has degree $\le m$ and at most one factor can have degree
 $<m$.

Suppose now that   $m$ is even. 
Then, since as remarked above all factors but at most one  have degree
 $m$, we should have $\mu_i=0$ for at least three of the factors,
 corresponding say  to $i=\alpha,\beta,\gamma$. But this would give a
 rational parameterization of the elliptic curve
 $y^2=(x-\xi_\alpha)(x-\xi_\beta)(x-\xi_\gamma)$, a contradiction.

Therefore $m$ is odd, which forces  $m=\deg a>2\deg q$, and in
 particular all of the said factors have the same degree $m$, so
 $\mu_i=1$ for all $i$ in \eqref{E.3}; note that this implies in
 particular that $\lambda-\xi_i$ divides $a(\lambda)- \lambda q^2(\lambda)$
 for each $i$.  

Now, since the degree of the whole product is  $2\deg b$, we must have
 $\deg (a(\lambda)-tq^2(\lambda))$ even, which implies $\deg a=1+2\deg
 q=1+2h$, say.
 It also follows that  $\deg c_i=h$. 

Let now $s_i:=(\lambda-\xi_i)\left({c_i(\lambda)\over q(\lambda)}\right)^2$.
  We have that
 $s_i-s_1=\xi_1-\xi_i$ is   constant; hence,
\begin{equation}\label{E.s}
s_i'=s_1',\qquad i=1,\ldots ,d.
\end{equation}
We compute 
$$
s_i'={c_i\over q^3}((c_i+2(\lambda-\xi_i)c_i')q-2(\lambda-\xi_i)c_iq')
={c_i\over
 q^3}\phi_i,
$$
say, where $\phi_i$ are non-zero polynomials of degree $\le 2h$
 (in fact $=2h$, as may be checked).  From \eqref{E.s} we find
 $c_i\phi_i=c_1\phi_1$. But the $c_i$ are pairwise coprime; hence
 $\prod_{i=2}^dc_i$ divides $\phi_1$. But $\phi_1$ is not zero, since
 $s_1$ is not constant, whence $(d-1)h\le \deg \phi_1\le 2h$, which
 implies $h=0$.  But then $q$ is constant and $\deg a=1$, giving a
 contradiction with the fact that all $\lambda-\xi_i$ divide
 $a(\lambda)-tq^2(\lambda)$. 
 This  concludes the proof.
\end{pf}

\section{Transcendence of sections and the logarithms}
\subsection{Elliptic schemes and the exponential map}\label{ell-sch}
Let $B$ be a smooth affine algebraic curve over $\C$,
 and $\pi:\, \cE\to B$ be an elliptic scheme.
 By this we mean that {\it each} fiber $\pi^{-1}\{t\}=\cE_t$ with
$t \in B$ is a smooth elliptic curve; in other words, the bad
reduction can arise only at the points at infinity of a completion of
$B$.

Every elliptic curve $\cE_t$, for $t\in B$, has a Lie algebra $\Lie(\cE_t)$,
which is a one-dimensional vector space. The  union of these lines
constitutes a line bundle $\Lie(\cE)\to B$ over $B$, which is
holomorphically trivial by the Oka-Principle
($B$ is a one-dimensional Stein manifold; cf., e.g., \cite{n16} Theorem 5.5.3).
 The exponential map
$\Lie(\cE_t)\to \cE_t$ has a kernel $\Lambda_t$, which is a discrete group
of rank two. These groups together define a  local system over $B$,
i.e. a sheaf  in abelian groups $\Lambda$, which is locally isomorphic to
the constant sheaf associated to the group $\Z^2$.

Recalling that $\cE_t$ too has the structure of an abelian group, so that
to the family $\cE\to B$ one can  associated the group-sheaf of its
holomorphic sections, we have the short exact sequence
\begin{equation}\label{exact1}
0\to \Lambda\to \Lie(\cE)\to \cE \to 0.
\end{equation}
From another perspective, we can view $\Lambda$ as a Riemann surface
covering $B$, $\Lie(\cE)$ as the total space of a line-bundle, i.e. an
algebraic surface fibered over $B$, and $\cE$ as an (open set of an)
elliptic surface.
Taking the long sequence in cohomology from \eqref{exact1},
 we obtain
\begin{equation}
0\to \Gamma(B,\Lambda)\to \Gamma(B,\Lie(\cE))\to \Gamma(B,\cE)\to
 \mathrm{H}^1(B,\Lambda)\to 
 \mathrm{H}^1(B,\Lie(\cE))=0. 
\end{equation}
The last zero is due again to the fact that $B$ is a one-dimensional
Stein manifold. Now, if the elliptic scheme $\cE\to B$ is not isotrivial,
 then no non-zero
period can be continuously defined in $B$; hence the term
$\Gamma(B,\Lambda)$ also vanishes. We finally get
\begin{equation}\label{seclog}
\frac{\Gamma(B,\cE)}{\exp(\Gamma(B,\Lie(\cE)))} \iso \mathrm{H}^1(B,\Lambda).
\end{equation}
The group $\Gamma(B,\cE)$ consists of the group of holomorphic sections:
it properly contains the Mordell-Weil group, formed by the rational
sections. The latter is a discrete group, since it injects into the
discrete group $\mathrm{H}^1(B,\Lambda)$ via the above projection
${\Gamma(B,\cE)}\to \mathrm{H}^1(B,\Lambda)$. In other words,
{\em no non-zero rational section of an elliptic scheme can admit a logarithm},
i.e. a lifting to $\Lie (\cE)$.

 We shall see that a holomorphic section of an abelian scheme
 $\sA \to B$ in general (and a semi-abelian scheme with an additional condition)
admitting a well-defined logarithm is transcendental or a constant section in its
$\C(B)/\C$-trace (see Theorems \ref{abl}, \ref{sabl}). 

{\bf N.B. } From the above discussion, it follows that the group of holomorphic
sections of an elliptic scheme is an extension of a finitely generated
group by an infinite dimensional vector space.

\subsection{Transcendency of sections}\label{trs}
Let $A$ be a semi-abelian variety  over $\C$ of dimension $n$; it is
the middle term in the exact sequence:
\begin{equation}\label{sa0}
0 \to \bG_m^l \to A \to A_0 \to 0, 
\end{equation}
where $A_0$ is an abelian variety.
 Let $\Lie (A) \to A$ be the Lie algebra of $A$, endowed with its
exponential map; analytically,
\begin{equation}\label{lie}
\Lie (A) \iso \C^n \to A=\C^n / \Gamma
\end{equation}
for a discrete subgroup $\Gamma$ (semi-lattice) of $\C^n$.

Let $B$ be a smooth affine algebraic curve.
 We consider the relative setting
of \eqref{sa0} over $B$:
\begin{equation}\label{sa1}
\begin{array}{ccccccccc}
0 & \to & \bG_{mB}^l & \to & \sA & \overset{\phi}{\lto} & \hskip-6pt \sA_0 & \to &
 0 \\
&  &  & \searrow & \,~ \downarrow\hskip-2pt\pi &
~\, \swarrow\hskip-2pt\pi_0 & & & \\
& & & & B & & & &  \\
\end{array}
\end{equation} 
Here we assume that $\pi_0: \sA_0 \to B$ is smooth without degeneration
and also  $d\pi$ is everywhere non-zero; 
in this case ,
 we say that $\pi: \sA \to B$ is {\em smooth}. 
After deleting possibly a finite number of points of $B$, we may reduce
the initial case to a smooth one.

As in \eqref{lie} we have the relative Lie algebra over $B$ and the
corresponding semi-abelian exponential 
\begin{equation}\label{exp1}
\varpi : \Lie(\sA) \lto \sA .
\end{equation}

For a point $t \in B$ we have
 $\sA_t=\pi^{-1} \{t\}  \iso \C^n/\Gamma_t$,
 where  $\Gamma_t$ is a semi-lattice.
Since, as above, the vector bundle $\Lie(\sA) \to B$ is analytically trivial
by Grauert's Oka-Principle ($B$ is a one-dimensional Stein manifold;
 cf., e.g., \cite{for} Theorem 5.3.1),
we can write
\begin{equation}\label{exp2}
\varpi : \Lie(\sA) \iso  B \times \C^n \ni (t, x)  \mapsto
 [(t,x)] \in \C^n/\Gamma_t=\sA_t
 \subset \sA.
\end{equation}

Let  $\sigma: B \to \sA$ be a holomorphic section of
$\pi:\sA \to B$.
If there is a lifting $\tilde\sigma: B \to \Lie(\sA)$ in \eqref{exp1}
with $\varpi \circ \tilde\sigma (t)=\sigma (t)~(t \in B)$, we call
$\tilde\sigma$ a {\em logarithm} of $\sigma$ (over $B$).

As for the case of elliptic schemes, already analyzed, logarithms do not
always exist (cf.\ \S\ref{ell-sch}).
 The semi-lattices $\Gamma_t$ define  a local system over $B$,
 and the existence of a logarithm for a section $\sigma$ is
obstructed by a cohomology class in the corresponding first cohomology
group of this local system (cf.\ \eqref{seclog}).

With reference to \eqref{sa1} we denote by $\sG_0$ the $\C(B)/\C$-trace of
$\sA_0$ with the quotient morphism $q_{0}: \sA_0 \to \sA_0/\sG_0$.
We have 
\begin{equation}\label{exct1}
 0 \to \bG_{mB}^l \to \sG_1 := \Ker \phi\circ q_0 \to \sA
~\overset{\phi}{\lto}~ \sA_0~
 \overset{q_0}{\lto}~ \sA_0/\sG_0 \to 0 , 
\end{equation}
and hence the exact sequence
\[
0 \to \bG_{mB}^l \to \sG_1  \to \sG_0 \to 0  .
\]
Thus, $\sG_1$ gives rise to a semi-abelian scheme over $B$.

Note that $\sG_0$ is defined over $\C$ and $\sG_0$ is isomorphic, as a
scheme over $B$, to a product $B \times G_0$ for an abelian variety
$G_0$;  $\sG_0\to B$ is the ``constant part'' of
the abelian scheme $\sA_0\to B$.
We say that a holomorphic section $\sigma:B \to \sA$ is
{\em $\sG_0$-valued constant} if $\phi\circ\sigma(B) \subset \sG_0\iso B \times G_0$
and $\phi\circ\sigma(t)=(t, x_0)$ with an element $x_0 \in G_0$.

We consider $\bG_{mB}^l$ in \eqref{sa1}.
Since the only complex affine algebraic model of $\bG_m$ is $\C^*$,
after a finite base change  we have
\begin{align}\label{triv}
 \bG_{m B}^l &\iso B \times (\C^*)^l,\\ \notag
\nonumber
0 \to B \times (\C^*)^l &\to \sG_1  {\lto} B \times G_0 \to 0
\quad (\hbox{over } B). 
\end{align}
We keep this reduction and the notation henceforth.

Taking a smooth equivariant toroidal compactification $T$ of $(\C^*)^l$,
 we have a fiber bundle
\begin{equation}\label{cpt}
 \bar\sA\, \underset{T}{\lto}\, \sA_0 ~(\to B).
\end{equation}
We then have the space $\Omega^1(\bar\sA, \log \del \sA)$
of logarithmic 1-forms with $\del\sA= \bar\sA\setminus \sA$
and
$\mathbf{T}(\bar\sA, \log \del \sA)$ of logarithmic vector
fields along the divisor $\del \sA$.

We consider the transcendency problem of a holomorphic section of
$\sA \to B$ with a logarithm. If $\sA \iso B \times A_1$
 (trivial family), then any constant section
 of $B \times A_1 \to B$ is rational and has a logarithm;
 this may happen in a subfamily
of $\sS \to B$ of $\sA \to B$ , even if $\sA \to B$ is non iso-trivial.

It is also to be noticed that a holomorphic section defined in a neighborhood
of a point of $\bar{B}\setminus B$ may locally have a 
non-constant logarithm there. But, globally we have:

\begin{lem}\label{logab}
Let $A_0$ be an abelian variety with an exponential map,
$\exp: \C^n \to A_0$.
Let $g: \Delta^*=\{0<|z| < 1\} \to A_0$ be a holomorphic map
with a logarithm $f: \Delta^* \to \C^n$ such that $g(z)=\exp f(z)$.
If $g(z)$ is holomorphically extendable at $0$ 
 as a map into $A_0$,
 then so is $f(z)$ as a vector-valued holomorphic function.

In particular, if $g:B \to A_0$ is a rational map with a logarithm,
 then $g$ is constant.
\end{lem}
\begin{pf}
Assume that $g:\Delta^* \to A_0$ is holomorphically extendable at $0$.
Then $f: \Delta^* \to \C^n$ is reduced to be bounded in a small
punctured neighborhood of $0$, and so Riemann's extension
implies that $f$ is holomorphically extendable at the puncture $0$.

Let $f$ be a logarithm of the rational section $g: B \to A_0$ and
 let $\bar B$ be a smooth compactification of $B$.
Since $g$ extends to a holomorphic map $\bar B\to A_0$,
 $f$ extends holomorphically over $\bar B$ as a vector-valued 
holomorphic function.
Hence, $f$ is constant and so is $g$.
\end{pf}

In view of Mordell-Weil over function fields (Lang-N\'eron) we
 have
\begin{thm}\label{abl}
Let $\sA_0 \to B$ be an abelian scheme and let $\sG_0$ be a $\C(B)/\C$-trace
of $\sA_0$. 
Let $\sigma: B \to \sA_0$ be a holomorphic section with a logarithm.
Then $\sigma$ is either $\sG_0$-valued constant or transcendental.
\end{thm}
\begin{pf}
Let $q: \sA_0 \to \sA_0(\C(B))/\sG_0$ be the quotient map.
Then $q \circ \sigma$ is a rational section of $\sA_0/\sG_0$ over $B$
with a logarithm.
By Lang-N\'eron, $\sA_0(\C(B))/\sG_0$ is finitely generated and hence
 discrete; in particular, no non-zero section of $\sA_0(\C(B))/\sG_0$ is
 infinitely divisible. Now, if a   rational section
 $\rho:B\to\sA_0(\C(B))/\sG_0$ admits a logarithm, this section is
 infinitely divisible in the holomorphic sense, i.e. for every integer
 $n$   there exists a holomorphic section $\rho_n$ with $n\cdot
 \rho_n=\rho$. This last equation is an algebraic one, so every solution
 is algebraic; since $\rho_n$ is well-defined on the whole of $B$, being
 algebraic it must be rational. It follows that $\rho$ is infinitely
 divisible in the Mordell-Weil group,
 and hence it is the $0$-section. 

Applying this fact to  $\rho=q\circ \sigma$ we obtain that $q\circ
 \sigma=0$ so   $\sigma: B \to \sG_0\iso B \times G_0$.
We write
\[
 \sigma(t)=(t, \exp f(t)),
\]
where $\exp: \Lie (G_0)\iso \C^{n_0} \to G_0$ is an
exponential map and $f: B \to \C^{n_0}$ is a vector-valued holomorphic
 function. By Lemma \ref{logab}, $f(t) \equiv a_0 \in \C^{n_0}$
 and $\sigma(t)=(t, x_0)$
with $x_0=\exp a_0$.
\end{pf}

To generalize the above results to semi-abelian varieties we need:
\begin{lem}\label{logtr}
Let $g(z)$ be a holomorphic function on a punctured disk
${\Delta}^*=\{z \in \C: 0 < |z| < 1\}$.
If $g(z)$ is not extendable at $z=0$ as a holomorphic function,
then $e^{g(z)}$ has an essential (isolated) singularity at $0$.
\end{lem}
 {\it Remark.}   In function theory, ``$e^{\mathrm{transcendental}}=
\mathrm{algebraic}$'' does not happen, while in numbers, $e^{\pi i}=-1$.
\begin{pf}
We distinguish two cases, according to the type of singularity of $g$ at $0$: 

(i) $g$ has a pole at $0$. Then in every punctured neighborhood of $0$,
 the real part $\re g(z)$ of $g(z)$ takes arbitrarily large
 positive numbers and arbitrarily small negative numbers,
so the function $e^{g(z)}$ tends to infinity on a sequence
 converging to $0$ and it also tends to $0$ on another such
 sequence. This can happen only if $e^{g(z)}$ has an essential
 singularity at $0$.

(ii) $g$ has an essential singularity at $0$. Then the image by $g$ of
 any punctured neighborhood of $0$ is dense, so again $g$ tends to two
 different values on sequences converging to $0$ (say it tends to $0$
 and to $1$) so $e^{g(z)}$ has two limits on different sequences. 
Thus, $e^{g(z)}$ has an essential (isolated) singularity at $0$.
\end{pf}

\begin{thm}\label{sabl}
Let $\pi: \sA \to B$ be a smooth semi-abelian scheme and let $\sG_1$ be as in
\eqref{exct1}. Assume that
$$
\sG_1 \iso B \times G_1
$$
with a semi-abelian variety $G_1$ over $\C$.
If a holomorphic section $\sigma: B \to \sA$ has a logarithm,
 then $\sigma$ is either transcendental or $\sG_1$-valued constant,
i.e., $\sigma(t)=(t, x_1)$ with an element $x_1 \in G_1$
through $\sG_1\iso B\times G_1$.
\end{thm}
\begin{pf}
By Theorem \ref{abl} and Lemmata \ref{logab}, \ref{logtr}.
\end{pf}

\section{Local smooth family}\label{lsf}
We would like to transpose the Nevanlinna theory for holomorphic
curves into semi-abelian varieties (cf.\ \cite{n81}, \cite{nw14} Chap.\ 6)
to a relative setting. 
\subsection{Jet space of holomorphic local sections}
Let $\Delta$ be the unit disk of the complex plane $\C$
with center $0 \in \C$.
 Let $t \in \Delta$ be the natural complex coordinate.
We consider a {\em smooth family} $\pi: \sA \to \Delta$ of semi-abelian
varieties of dimension $n$ with its zero section:
 $  \Delta \ni t \mapsto 0_t \in \sA_t=\pi^{-1}\{t\}$, $t \in \Delta$.

Let $\bar\sA$ be a relative toroidal compactification of $\sA$
(cf.\ \eqref{cpt}).
Let $J_k(\bar\sA, \log \del \sA)$ denote the $k$\,th logarithmic
jet space over $\bar\sA$ along $\del\sA$, and let
$$
\pi_k: J_k(\bar\sA, \log \del \sA) \to \bar\sA
$$
be the natural projection.

Let $J_{k}(\sA/\Delta) (\subset J_k(\sA))$
 denote the space of $k$-th jets
of holomorphic local sections $f$ of $\pi: \sA \to \Delta$ such
that $\pi\circ f(t)=t$. 

In \eqref{exp2} we write $x=(x_1, \ldots, x_n)$ with the natural complex
coordinates. Then, $\eta_j:=d x_j$ ($1 \leq j \leq n$) give rise to
  elements of the space $\Omega^1(\bar\sA, \log\del\sA)$
 of logarithmic $1$-forms and 
\begin{equation}\label{hframe}
\{dt, \eta_1, \ldots, \eta_n\}
\end{equation}
forms the frame over $\bar\sA$.

For a jet element $j_k(f) \in J_k(\sA/\Delta)_{f(t)}$ ($t \in \Delta$) we
set
\begin{equation}\label{der}
 f^* \eta_j=f'_j dt , \quad 1 \leq j \leq n.
\end{equation}
Then we have
\begin{align}\label{jets}
j_1(f)(t) &=(f(t); 1, f'_1(t), \ldots, f'_n(t)), \\ \nonumber
j_2(f)(t) &=(j_1(f)(t); 0, f''_1(t), \ldots, f''_n(t)), \\ \nonumber
& ~ \vdots \\ \nonumber
j_k(f)(t) &=(j_{k-1}(f)(t); 0, f^{(k)}_1(t), \ldots, f^{(k)}_n(t)).
\end{align}
In this way we have the trivializations
\begin{equation}\label{prod}
 J_k(\sA/\Delta) \iso \sA \times \{(1,0, \ldots, 0)\} \times \C^{nk}
\iso \sA \times \C^{nk}.
\end{equation}
Let
\[
  I_k:   J_k(\sA/\Delta) \to  \C^{nk}
\]
be the jet projection, which extends holomorphically to
the relative logarithmic jet space\break
 $J_k(\bar\sA/\Delta,\log \del \sA)$.
We set the relative jet projection with respect to the frame \eqref{hframe}
\begin{equation}\label{jetp}
 \tilde{I}_k=(\pi_k, I_k):  J_k(\bar\sA/\Delta, \log \del \sA) \to 
 \Delta \times \C^{nk}.
\end{equation}
Note that $ \tilde{I}_k$ is proper.

\subsection{A relative exponential map}
We keep the notation above. We consider the abelian integration
\begin{equation}\label{abint}
x_t \in \sA_t \to \left( \int_{0_t}^{x_t} \eta_1, \ldots,
 \int_{0_t}^{x_t} \eta_n \right) \in \C^n.
\end{equation}
We denote by $\Gamma_t$ the semi-lattice generated by the periods
of \eqref{abint}.
We then have a relative exponential map 
\begin{align*}
 \exp_\Delta : \Lie (\sA) \iso   \Delta \times \C^n \ni (t, x) & \mapsto
(t, [x]) \in \{t\} \times \C^n/\Gamma_t =\sA_t \subset \sA, \\
  \pi\circ \exp_\Delta(t, x) &=t.
\end{align*}

For an element $w \in \C^n$ we have an {\em action} ``$w\cdot$''
associated with \eqref{hframe} by
\begin{equation}\label{act}
w\cdot : (t, [x]) \in \{t\} \times \C^n/\Gamma_t =\sA_t \to
(t, [x+w]) \in \{t\} \times \C^n/\Gamma_t =\sA_t.
\end{equation}  

\subsection{Nevanlinna theory of holomorphic sections over the punctured disk}
The notation is kept. We follow \cite{n81}.

Let $\Delta^*=\Delta \setminus \{0\}$
be the punctured disk. We consider a holomorphic section
\[
 f: \Delta^* \lto \sA,\qquad \pi \circ f (t)=t, ~ t \in \Delta^*.
\]

Let now $\omega$ be a real
$(1,1)$-form  on $\bar\sA$ and
 let $r_0>1$ be any fixed real number.
For $r>r_0$ we define the {\it order function of $f$ with respect to $\omega$}
 by
\begin{equation}\label{ordf}
 T_f(r;  \omega)=\int_{r_0}^{r}
 \frac{ds}{s} \int_{\{1/s < |t| <1/r_0\}} f^*\omega, \quad r > r_0.
\end{equation}
Let $\omega_0$ be a hermitian metric form on $\bar\sA$.
Then there is a  constant $C>1$ such that
\begin{equation}\label{orderf}
 C^{-1} T_f(r; \omega_0) \leq  T_f(r; \omega) \leq  C T_f(r;  \omega_0),
\quad r > r_0.
\end{equation}
The above $C$ depends on the choice of $r_0>0$ in general.

It is noted:
\begin{prop}\label{Ologr}
Let $\omega_0$ be a hermitian metric form on $\bar\sA$.
A holomorphic section $f: \Delta^* \to \sA$
 is holomorphically extendable at $0$ 
as a map into $\bar\sA$ if and only if
\begin{equation}
 \lowlim_{r \to \infty} \frac{T_f(r; \omega_0)}{\log r} < \infty.
\end{equation}
\end{prop}

Let $\sD$ be a relative effective divisor on $\sA/\Delta$
which is extendable to a divisor $\bar\sD$ on $\bar\sA/\Delta$.
We call such $\sD$ a {\em relative algebraic} divisor on $\sA/\Delta$. 
Let $\sL=\sL(\bar\sD)$ denote the line bundle over
$\bar\sA/\Delta$ determined by $\bar\sD$ with
a section $\sigma$ such that the divisor $(\sigma)$ defined by
 $\sigma$ satisfies $(\sigma)=\bar\sD$. Let $\| \cdot \|$ be a
hermitian metric in $\sL$, and let $\omega_\sL$ be the Chern curvature form
of the hermitian metric.
For a holomorphic section $f: \Delta^* \to \sA$ with
$f(\Delta^*) \not\subset \supp\,\sD$, we define the counting functions
of the pull-back divisor $f^*\sD$ by
\begin{align*}
n(s, f^* \sD) &= \sum_{1/s <|\zeta|< 1/r_0} \deg_\zeta f^* \sD, \quad
 s>r_0,\\
 N(r, f^* \sD) &=\int_{r_0}^{r} \frac{n(s, f^* \sD)}s ds, \quad r> r_0.
\end{align*}
Replacing $\deg_\zeta f^*\sD$ above by $\min\{\deg_\zeta f^*\sD, k\}$
($k \in \N$),
we have the corresponding (truncated) counting functions denoted by
\[
 n_k(s, f^*\sD), \quad  N_k(r, f^*\sD).
\]

We set $\Gamma(r)=\{t=1/(re^{i\theta}): 0 \leq \theta\leq 2\pi\}$
parameterized by $\theta$, and the {\em proximity function}
\[
 m_f(r, \sD)=\int_{\Gamma(r)} \log \frac1{\|\sigma\circ f \|}\, \frac{d\theta}{2\pi}.
\]
We have:
\begin{thm}[\rm First Main Theorem (cf.\ \cite{n81} (1.4))]\label{fmt}
Let the notation be as above. Then
\begin{align}\label{fmteq}
 T_f(r ; \omega_\sL) &=N(r, f^*\sD)+m_f(r, \sD) -m_f(r_0, \sD)+
(\log r)\int_{\Gamma(r_0)} d^c \log \|\sigma\circ f \|^2 \\
\nonumber
& =N(r, f^* \sD)+m_f(r, \sD) + O(\log r), \quad r>r_0,
\end{align}
where $d^c=(i/4\pi)(\delbar - \del)$.
\end{thm}

\begin{rmk}\label{ordfbig}\rm
\begin{enumerate}
\item
When $f(\Gamma(r_0))\cap \supp\, \sD \not= \emptyset$, the last term
of \eqref{fmteq}, $\int_{\Gamma(r_0)} d^c \log \|\sigma\circ f \|^2$
should be taken as a principal-value integration. We may also take
$r_0>1$ so that $f(\Gamma(r_0))\cap \supp\, \sD = \emptyset$.
Then the integrand is smooth on $\Gamma(r_0)$.
\item
(Cf.\ \eqref{orderf})
If $\sL=\sL(\bar\sD)$ is relatively big on $\bar\sA/\Delta$, then there is
a positive constant $C$ such that
\[
 C^{-1}T_f(r; \omega_\sL)+O(\log r) < T_f(r; \omega_0)
<  CT_f(r ;  \omega_\sL)+O(\log r).
\]
\end{enumerate}
\end{rmk}

Let $f: \Delta^* \to \sA$ be a holomorphic section
and the frame \eqref{hframe} be given.
Recall we have defined the first derivatives $f'_j(t)$ of $f$ by \eqref{der},
and hence the $k$-th ($k \in \N$, positive integers)
 derivatives $f^{(k)}_j (t)$,
which are holomorphic functions on $\Delta^*$.
We then have the ``lemma on logarithmic derivatives'':
\begin{lem}[\cite{n77}, \cite{n81}]
Let the notation be as above. Then we have
\[
 \int_{\Gamma(r)} \log^+ \left| f^{(k)}_j \right| 
~ \frac{d\theta}{2 \pi}=S_f(r ; \omega_0),
\quad r >r_0,\quad k \in \N ,
\]
where $S_f(r ; \omega_0)=O(\log^+ T_f(r ; \omega_0))+O(\log r)||$, called
a small term in Nevanlinna theory.
\end{lem}

\subsection{Relative second main theorem}
Let $\sA \to \Delta$ and  $\bar\sA \to \Delta$ be a smooth family of
semi-abelian varieties and its relative toroidal compactification
 as above.
We consider a relative algebraic reduced divisor $\sD$ on $\sA/\Delta$;
i.e., there is a relative reduced divisor $\bar\sD$ on $\bar\sA$ such
that $\sD=\bar\sD \cap \sA$.
We identify it with its support.
For a given holomorphic section $f: \Delta^* \to \sA$
we deal with a problem to obtain a ``Second Main Theorem''
with respect to  $\sD$.

We refer to the {\em relative Zariski} topology on $\bar\sA$
in the sense that closed subsets $Z ~(\subset \bar\sA)$ are
 analytic subsets of $\bar\sA$;
hence, the fibers $Z_t ~(t \in \Delta)$ are algebraic
subsets of $\bar\sA_t$. Then, it induces the {\em relative Zariski} topology
on the open subset $\sA$ of $\bar\sA$. It is noticed that
a relative Zariski-closed subset $Y$ of $\sA$ is not merely an
analytic subset of $\sA$ but the fibers $Y_t ~(t \in \Delta)$
are algebraic subsets of $\sA_t$.
We define similarly the relative Zariski topology on the jet space
$J_k(\bar\sA/\Delta, \log \del\sA)$ and its open subset
$J_k(\sA/\Delta)$:
Here one notes that the restrictions
$J_k(\bar\sA/\Delta, \log \del\sA)|_{\bar\sA_t} ~(t \in \Delta)$
are  affine fiber bundles over $\bar\sA_t$ with some $\C^r$ as fibers,
 which is compactified by $\P^r(\C)$.

Let $f: \Delta^* \to \sA$ be a holomorphic section.
We denote by
 $X_k(f)~ (\subset J_k(\bar\sA/\Delta, \log \del\sA), k \geq 0)$
the relative Zariski-closure of the image of the $k$-th jet lift
$J_k(f): \Delta^* \to J_k(\sA/\Delta)\hookrightarrow J_k(\bar\sA/\Delta, \log \del\sA)$
of $f$.
 For $k=0$ we set
\begin{equation}\label{xf}
X(f)=X_0(f) ~\subset \bar\sA . 
\end{equation}
\begin{defn}\rm We say that
$f: \Delta^* \to \sA$ is {\em non-degenerate}
if $X(f)=\bar\sA$; otherwise, $f$ is {\em degenerate}.
\end{defn}

If $f: \Delta^* \to \sA$ is  non-degenerate, $f$ 
is  not extendable at $0$ as a map into $\bar\sA$.

To prepare some technical lemmata, we need to fix the frame \eqref{hframe}.

\begin{lem}\label{jetdiff}
Assume that $f: \Delta^* \to \sA$ is non-degenerate.
Then for all sufficiently large $k \in \N$  
\[
 \tilde{I}_k(X_k(f)) \cap \tilde{I}_k(\sD/\Delta)\not=
 \tilde{I}_k(X_k(f)).
\]
\end{lem}

\begin{pf} Cf.\ \cite{n81}, \cite{nw14} Lemma~6.3.2.
\end{pf}

As in the proof of \cite{nw14} p.~227, we have:
\begin{lem}\label{lemder}
Assume that $f: \Delta^* \to \sA$ is 
non-degenerate.
Then we have
\[
 m_f(r ; \sD)=S_f(r ; \omega_0), \quad r>r_0.
\]
\end{lem}

Combining this with Theorem \ref{fmt} we obtain
\begin{thm}[\rm Second Main Theorem]\label{smt}
Let $\sD$ be a relatively algebraic big reduced divisor on $\sA/\Delta$
 and let $f: \Delta^* \to \sA$ be a non-degenerate holomorphic section.

 Then there is a relative compactification $\bar\sA/\Delta$ of
 $\sA/\Delta$ together with $\sL=\sL(\bar\sD)$,
independent of $f$, 
and a natural number $k \in \N$ such that
\[
 T_f(r ; \omega_\sL)=N_k(r, f^*\sD)+S_f(r ; \omega_\sL), \quad r>r_0.
\]
\end{thm}
\begin{pf}
Cf.\ \cite{nw14} \S6.4.
\end{pf}

\begin{rmk}\rm
\begin{enumerate}
\item
Note that the order function $T_f(r ; \omega_\sL)$ depends on the choice
of the relative compactification $\bar\sA/\Delta$.
\item
The requirement on $\bar\sA$ above is such that for a general point
$t \in \Delta$,
$\sA_t$ acts on $\bar\sA_t$ and $\bar\sD_t$ contains no fixed point
(\cite{nw14}, Corollary 5.6.7); this is an open property in the
     parameter $t$.
\end{enumerate}
\end{rmk}

\begin{thm}[Big Picard]\label{bp}
Let $f: \Delta^* \to \sA$ and $\sD$ be as in Theorem \ref{smt}.
Then, $f$ intersects $\sD$ infinitely many times in an arbitrarily
 small (punctured) neighborhood about $0$.
\end{thm}
\begin{pf}
Suppose that $f(\Delta^*) \cap \sD$ is finite.
Since  $\sD$ is a relatively algebraic big reduced divisor on $\sA/\Delta$,
it follows from Theorem~\ref{smt} that
\[
 T_f(r; \omega_0)=S_f(r ; \omega_0), \quad r>r_0.
\]
Then, it is immediate that $T_f(r, \omega_0)=O(\log r)||$.
By Proposition~\ref{Ologr} we conclude that $f$ is extendable
at $0$ as a map into $\bar\sA$, and hence $f$ cannot be non-degenerate.
\end{pf}

\begin{rmk}\rm
By Corvaja-Noguchi \cite{cn} Theorem 5.2,  a non-degenerate
entire curve on a semi-abelian variety intersects an ample reduced divisor $D$
so that the points of the intersections {\em on} $D$ is Zariski-dense in $D$.
Together with the smoothness assumption of $\sA \to \Delta$ at $0$,
 it is interesting to ask:
\begin{quest}\rm
Is $f(\Delta^*) \cap \sD$ of Theorem \ref{bp} 
relative Zariski-dense in $\sD$?
\end{quest}

\begin{quest}\label{deg1}\rm
 Is it possible to allow the above $\sA \to \Delta$ of 
 Theorems \ref{smt} and \ref{bp} to degenerate at $0 \in \Delta$? ~(Cf.\
Problem \ref{deg2}.)
\end{quest}
\end{rmk}

\subsection{Remarks to the Zariski-local case }\label{rmkalg}
Let $\bar{B}$ be a smooth projective algebraic curve.
Let $\pi:\sA \to \bar{B}$ and $\pi:\bar\sA \to \bar{B}$ be a smooth family of
$n$-dimensional semi-abelian varieties over $\bar{B}$ and its relative toroidal
compactification as in \S\ref{trs}.

We consider a transcendental holomorphic section 
$$
f: B  \lto \sA,
$$
where $B$ is an affine open set of $\bar B$. Set $S=\bar B \setminus B$.
If $f$ is transcendental, then there is a point of $S$ at which $f$
 is {\em not} extendable as a map into $\bar\sA$.

To deal with the local Nevanlinna theory for such $f$
 as above obtained in the previous subsections,
we introduce a rational function $\tau$ on $\bar B$
such that $\tau$ is holomorphic on a Zariski-open neighborhood $B'$ of
$S$ in $\bar B$ and has a zero of order $1$ at every point of $S$.
Let $B'(r'_0) ~(r'_0>0)$ denote the union of connected components of
$\{x \in B': |\tau(x)|<1/r'_0\}$ containing the points of $S$.
Taking and fixing a large $r'_0 $, we have
\[
 d \tau(x) \not =0, \qquad x \in B'(r'_0).
\]
For simplicity, we denote by $B$ the Zariski-open set $B \cap B'$, and set
\[
 B(r) =\{x \in B: 1/r <|\tau(x)|<1/r_0 \},
\]
where $r_0>r'_0$ is any fixed number.

We use $\tau$ for the parameter $t$ of \eqref{hframe}.
For $\eta_j ~(1 \leq j \leq n)$ of \eqref{hframe} we define also them
by rational differentials on $\bar\sA$ with logarithmic poles
on $\del\sA=\bar\sA\setminus \sA$,
so that
\[
 \{d\tau, \eta_1, \ldots, \eta_n\}
\]
forms a holomorphic frame over $\bar\sA|_{B'(r'_0)}$, where $r'_0$
is replaced by a larger one if necessary.

We then define the order functions, counting functions
 and proximity functions as in the former subsections;
e.g.,  the order function $T_f(r; \omega)$ of $f$ with respect
to a hermitian metric form $\omega$ on $\bar\sA$ is defined by
\[
 T_f(r; \omega)=\int_{r_0}^r \frac{ds}{s}
\int_{B(s)} f^*\omega, \quad r>r_0.
\]

Replacing the relative Zariski topology by the (ordinary) Zariski
topology, we apply the arguments given in the previous subsections
for $f: B \to \sA$. We then obtain in particular
 the {\em Second Main Theorem \ref{smt} for $f: B \to \sA$}.

\section{About Problem \ref{st1}}
We deal with the semi-abelian case.
Let $\bar B$ be a smooth projective algebraic curve.
Here we need a rather strong assumption such that
$\sA \to \bar B$ is smooth. 

\begin{thm}\label{smooth}
Let $\sA \to \bar B$ be a smooth family of semi-abelian varieties over
$\bar B$,  and let $\sD$ be a relatively  big
reduced divisor on a relative compactification $\bar\sA/\bar B$.
 Let $B$ be any non-empty Zariski-open subset of $\bar B$, and
let $f: B \to \sA$ be a holomorphic section.
If $f(B') \cap \sD$ is finite,
 then $f$ is degenerate. 
\end{thm}

\begin{pf}
By the Second Main Theorem \ref{smt} and \S\ref{rmkalg}.
\end{pf}

\begin{rmk}\rm
(i) We insist that in the above example, the fibration is smooth on the
 projective curve, i.e. without bad reduction. Abelian schemes of that
 kind can be obtained as jacobians of so-called Kodaira fibrations (see
 \cite{cat}); examples with semi-abelian schemes arise e.g. from
 generalized jacobians of logarithmic Kodaira fibrations (definition 5
 in \cite{cat}).

(ii) In the theorem above, the family $\sA \to B$ may degenerate
at a point where $f$ is already defined holomorphically as a map into $\bar\sA$.
What is essentially excluded in the theorem
 is the case when the family $\sA \to B$ degenerates
at a possibly (essentially) singular point of $f$.
\end{rmk}

\section{About Problem \ref{st2} and \ref{st3}}
\label{prob1-2}
\subsection{On Problem \ref{st2}}\label{transc}

{\bf (a)} Example. We first give the example mentioned after Problem \ref{st2}.
Let $C$ be a smooth projective curve of genus $g \geq 2$,
and let $A$ be its Jacobian variety with an embedding $\eta: C \to A$.
Let $\exp_A: \C^2 \to A$ be an exponential map. Take an affine
line $L \subset \C^2$ such that the image $\exp (L)$ is Zariski-dense
in $A$. Let $C' \subset C$ be an affine open subset with a non-constant
rational holomorphic map
 $\varphi: C' \to \C \iso L$. 
Let us set $\sA:=C' \times A^2$
 with the first projection 
$\sA \to C'$ to the base space $C'$. With $\psi(x)=\exp_A (\varphi(x))$
we have a section
\begin{equation}\label{transc1}
 \sigma: x \in C' \lto (x, \eta(x), \psi(x)) \in \sA.
\end{equation}
By definition (cf. \cite{n18}), a non-constant holomorphic map
 $f$ from $C'$ into $A$ is {\em strictly transcendental} if
for every abelian subvariety $A_1$ of $A$, the composed map
$q_{A_1} \circ f: C' \to A/A_1$ with the quotient map
 $q_{A_1}: A \to A/A_1$ is {\em either transcendental or constant}.

Then, $\psi: C' \to A$ is strictly transcendental
 due to \cite{n18}.
Let $G(\eta) \subset C \times A$ denote the graph of $\eta$
and set $G'(\eta)=G(\eta) \cap (C' \times A)$. Then,
the Zariski-closure $X(\sigma)$ (cf.\, \eqref{xf}) of
 the image of $\sigma$ is given by
\[
 X(\sigma)=G'(\eta) \times A.
\]
Thus, $X(\sigma)$ contains a
translate of a subgroup, $\{w\} \times A$ with $w \in G'(\eta)$, 
 however it is not itself a translate of a subgroup.

\smallskip
{\bf (b)}  We use the notation defined in \S\ref{trs}.
Let $\pi: \sA \to B$ be an algebraic family of
 semi-abelian varieties over an affine curve $B$
   and let $\bar\sA \to B$ be the relative compactification over $B$.
 For a Zariski-closed subset $X \subset \sA$ we consider the 
 the set-theoretic stabilizer group of $X$, which is
\begin{equation}
 \bigcup_{t \in B}\{a \in \sA_t: a+X_t=
 X_t\} \quad (X_t=X \cap \sA_t).
\end{equation}
It is a closed subset of $\sA$; each fiber over any point $t\in B$ is an
algebraic subgroup of the fiber $\sA_t$; its dimension is upper
semi-continuous, i.e. either constant or admitting jumps at a finite set
of points. Removing these points, we let $St(X)$ be the Zariski-closure
of the set-theoretic stabilizer outside this exceptional set. It is the
total space of a group-scheme over $B$, so that we have the quotients,
\begin{equation}
X/ St(X) \subset \sA/St(X).
\end{equation}

We set the relative dimension of $St(X)$ by
\[
 \dim_B St(X)= \dim St(X) -\dim B=\dim St(X) -1.
\]

For a holomorphic section $f:B  \to \sA$ we set (cf.\ \eqref{xf} for $X(f)$)
\begin{equation}\label{stb}
 St(X(f))=St(X(f) \cap \sA), \qquad 
 \dim_B St(X(f))= \dim_B  St(X(f)\cap \sA).
\end{equation}

\begin{thm}\label{smooth2}
Let $\pi: \sA \to \bar B$ be a smooth family of semi-abelian varieties
over  $\bar B$.
Let  $f: B \to \sA$ be a transcendental  holomorphic section over
 an affine open subset $B \subset \bar B$.
 Then  $ \dim_B St(X(f))>0$.
\end{thm}
\begin{pf}
We take the restriction of
 the jet projection $\tilde{I}_k$ defined by \eqref{jetp} to
 $X_k(f) (\subset J_k(\bar\sA/B,$ $\log \del \sA))$
\[
 \tilde{I}_k|_{X_k(f)}: X_k(f) \lto \Delta \times \C^{nk}.
\]
Since $\tilde{I}_k$ is proper, the image $Y_k=\tilde{I}_k(X_k(f))$
is an analytic subset of $\Delta \times \C^{nk}$ which is algebraic
in $\C^{nk}$-factor. We fix arbitrarily a reference point $t_0 \in
 \Delta^*$ such that $j_k(f)(t_0)$ is a non-singular point of $X_k(f)$,
and consider the differential between the holomorphic tangent spaces
\[
 d\tilde{I}_k|_{X_k(f)}: \mathbf{T}(X_k(f))_{j_k(f)(t_0)} \to
\mathbf{T}(\Delta \times \C^{nk})_{\tilde{I}_k(j_k(f)(t_0))}\iso
 \C \times \C^{nk}.
\]
By \eqref{prod} we have
\[
 \mathbf{T}(X_k(f))_{j_k(f)(t_0)} \subset \mathbf{T}(\Delta)_{t_0}\times
(\mathbf{T}(\sA_{t_0}) \times \mathbf{T}(\C^{nk}))_{j_k(f)(t_0)} \iso
\C \times \mathbf{T}(\sA_{t_0}) \times \C^{nk}.
\]
With this product structure the kernel $\Ker d\tilde{I}_k|_{X_k(f)}$
satisfies
\[
W_k:= \Ker~ d\tilde{I}_k|_{X_k(f)} \subset \mathbf{T}(\sA_{t_0})_{f(t_0)} \iso \C^n,
\]
where the other tangent components are $0$.
Since $W_{k} \supset W_{k+1}$, they stabilize at some
 $W_{k_0}=W_{k_0+1}=\cdots$.
If $W_{k_0}=\{0\}$, we then deduce in the same way as in \cite{nw14} \S6.2 that
$$
T_f(r; \omega_0)=S_f(r ; \omega_0),
$$
so that $f$ is extendable over $\bar B$ as a map into $\bar\sA$; i.e.,
 $f$ is {\em not} transcendental.
Hence $W_{k_0} \not= \{0\}$. Since a tangent vector $v \in W_{k_0}$
is tangent to $X(f)$ at $f(t_0)$ with infinite order, we see that
$v (\in \Lie(\sA/B))$ is tangent to $X(f)\cap \sA$
 at all points of $X(f) \cap \sA$ (cf.\ \cite{nw14} \S6.2).
Thus, $\dim_B St(X(f))>0$.
\end{pf}

\smallskip
\begin{cor}\label{picext}
Let $f: B \to \sA$ be as in Theorem \ref{smooth}
and let $q: \sA \to \sA/St(X(f))$ be the quotient map.
Then the composite $g:=q\circ f: B \to \sA/St(X(f))$ is rational.
\end{cor}

\subsection{On Problem \ref{st3} }
Now, let $\pi: \sA \to \bar B$ be an abelian scheme defined over
a smooth projective algebraic curve $\bar B$.
Let $B \subset \bar B$ be an affine open subset.
We consider a transcendental holomorphic section $f: B \to \sA$.
Let $X(f)$ be the Zariski closure of $f(B)$ in $\sA$ and let
$St(X(f))$ be the stabilizer of $X(f)$  defined by \eqref{stb}.
Let $\sG$ be an abelian subscheme of $\sA$ over $\bar B$ and
let $q_{\sG}: \sA \to \sA/\sG$ be the quotient morphism.
We consider the induced composite map $q_{\sG} \circ f: B \to \sA/\sG$.
Then, as in Theorem \ref{abl} it makes sense to say whether
$q_{\sG} \circ f$ is a
 $\C(B)/\C$-trace (of $\sA/\sG$) valued constant section or  not.
Following to the notion of ``strict transcendency'' defined in
\cite{n18}, we give:

\begin{defn}\label{strtr}\rm
We say that a transcendental holomorphic section $f: B \to \sA$ is
{\em strictly transcendental} if for every abelian subscheme $\sG$
of $\sA$ the induced holomorphic section $q_{\sG} \circ f: B \to
 \sA/\sG$
is either transcendental or a $\C(B)/\C$-trace valued constant section.
\end{defn}

\begin{thm}\label{gcoset}
Assume that $\sA \to \bar B$ is smooth.
Let $f :B \to \sA$ be a  strictly transcendental holomorphic section.
Then, after a finite base change, $X(f)$ is a translate of an abelian
subscheme $\sG$ of $\sA$ by a rational section $\sigma: B \to X(f)$
such that $q_\sG \circ f: B \to \sA/\sG$ is a $\C(B)/\C$-trace valued
constant section.
\end{thm}

\begin{pf}
We set $X=X(f)$ and $\sG=St(X)$.
By the assumption and Corollary \ref{picext}
\begin{equation}\label{constsec}
q_\sG \circ f: B \lto \sA/\sG
\end{equation}
is a $\C(B)/B$-trace valued constant section.
Thus for general $t \in B$ except for finitely many, the fibers satisfy
\[
 X_t= f(t)+ \sG_t.
\]
We take a finite base change $\tilde B \to B$ so that
there is a rational section $\tilde\sigma: \tilde B \to \tilde X$
with the lift $\tilde X$ of $X$. We denote by $\tilde\sA$
(resp.\ $\widetilde\sG$\,) the lift of $\sA$ (resp.\ $\sG$).  Thus, we have
\[
 \tilde X= \tilde\sigma + \widetilde\sG.
\]
It follows that
 $q_{\tilde\sG} \circ \tilde\sigma: \tilde B \to
 \tilde\sA/\widetilde\sG$
is the lift of the section \eqref{constsec} and hence a
 $\C(\tilde B)/\C$-trace valued constant section.
\end{pf}

We would like to stress  that we are assuming the smoothness of
the family $\pi: \sA \to \bar B$ over $\bar B$.
 It is interesting to ask:
\begin{quest}\rm
(i) Is it sufficient in Theorem \ref{gcoset} to assume
 the smoothness condition for $\sA$ only
 over the affine open $B$?

(ii) How do we deal with the semi-abelian case?
\end{quest}

\section{Problems and remarks}
{\bf (a)} As mentioned in \S\ref{intr}, the Zariski-closure
of the image of an entire curve $f: \C \to A$ in an abelian variety $A$
is a translate of an algebraic subgroup.
 It is still unknown what should be the topological closure
 in the sense of the complex topology, of the image
$f(\C)$. One might ask the following:

\begin{quest}\label{top}\rm
Let $f:\C\to A$ be an entire curve in a complex abelian variety. Is it
 true that the {\em topological closure} of $f(\C)$ is a translate of a real
 Lie subgroup of $A(\C)$?
\end{quest}

{\bf (b)} Related to Question \ref{deg1}, one may raise:

\begin{probl}\label{deg2}\rm
Extend the results of \S\ref{transc}
to the case where $\pi: \sA \to \bar B$ degenerates at finitely many points
of $\bar B$.
\end{probl}

\section{Heights--observations}

Let $B$ be an affine algebraic smooth curve.
For the heights of an abelian or elliptic scheme $\cE \to B$,
 we have to consider  three different settings:
\begin{enumerate}
\item
The height of a family  $\cE \to B$.
\item 
The height of a fiber $\cE_b$ for each $b \in B$.
\item
The height of a section $\sigma: B \to \cE$ of $\cE \to B$.
\end{enumerate}

{\bf (i) Height of elliptic schemes. } 
Let $\cE \to B$ be an elliptic scheme with a smooth projective algebraic
compactification $\bar B$ of $B$. The moduli space of elliptic curves
is an orbifold complex space $\bH/\Gamma(1)$ with
$\lambda:\bH=\{z \in \C : \im z>0\} \to \bH/\Gamma(1)$.
 We have a holomorphic map
$\phi_{\cE}: B \to \bH/\Gamma(1)$, which has a local lift
$\tilde\phi_{\cE U}:U \to \bH$
from a neighborhood $U$ of every point of $B$.
Since the Poincar\'e metric $\omega_{\rP}$ (or the Bergman metric in the case
of Siegel space) is invariant by $\Gamma(1)$, we have the pull-back
$\phi_{\cE}^*\omega_\rP=\tilde\phi_{\cE U}^*\omega_\rP$.
We define the {\em height} of $\cE \to B$ by
\begin{equation}
 h(\cE/B)=\int_B \phi_{\cE}^*\omega_\rP \in \R_{\geq 0}.
\end{equation}

Note that the integral is finite, because $\phi_{\cE}^*\omega_\rP$
has at most the Poincar\'e growth at every point of $\bar B \setminus B$.
 Also, we have
$h(\cE/B)=0$ if and only if the scheme $\cE/B$ is isotrivial. 
After the normalization of the curvature of $\omega_\rP$ we have
(by Schwarz or by the comparison of Kobayashi hyperbolic metrices)
\begin{equation}
 h(\cE/B) \leq 2g(B)-2 + \#(\bar{B}\setminus B) =\chi(B),
\end{equation}
where $\#(\bar{B}\setminus B)$ denotes the cardinality, 
 so  $\chi(B)$ is the Euler characteristic of the affine curve $B$.
The same holds for a family of principally polarized abelian varieties, up to replacing
the upper half plane by the Siegel space. 

{\bf (ii) Height of $\cE_b$. } Consider a non-isotrivial elliptic scheme
$\cE\to B$; identifying the Lie algebra 
$\Lie(\cE)$ with $B\times \C$, we can identify the period
lattice $\Lambda_b$ over each point $b$ with a
 lattice  in $\C$, and calculate its volume $V(b)$.
 Then we can define the {\em height} of $\cE_b$ to be
$$
h(\cE_b) = \frac1 {V(b)}.
$$ 
This function on $B$ depends on the trivialization of the line bundle
$\Lie(\cE)$. A canonical choice can be done by considering first
the Legendre scheme, whose base is $B_0=\P^1\setminus\{0,1,\infty\}$.
There one disposes of a canonical choice for the fundamental periods
$\rho_1,\rho_2$, given by Gauss' hypergeometric series defined in the
region of $B_0$ where $|\lambda|<1$ and $|1-\lambda|<1$:
\begin{equation}\label{E.Gauss}
\rho_1(\lambda)=\pi \sum_{n\geq 0} {{1/2}\choose{n}}^2 \lambda^n,\qquad
 \rho_2(\lambda)= i\pi \sum_{n\geq 0} {{1/2}\choose{n}}^2(1-\lambda)^n.
\end{equation} 
Then the volume of the lattice equals
\begin{equation*}
V(\lambda)=\frac{i}{2}\det \left(\begin{matrix} \rho_1 &
 \overline{\rho_1}\\ \rho_2 &
 \overline{\rho_2}\end{matrix}\right)=
\frac{i|\rho_1|^2}{2}(\bar{\tau}-\tau)=\frac{\Im(\tau)\cdot |\rho_1|^2}{2}.
\end{equation*}
Here $\tau=\tau(\lambda)=\rho_2/\rho_1$.
 Recall that the $\lambda$ function is expressed in terms of $\tau$ as 
 \begin{equation*}
 \lambda = 16\, q - 128\, q^2 +\ldots,\qquad \mathrm{for}\ q=e^{\pi i\tau}.
 \end{equation*}
Now, given any non-isotrivial elliptic scheme, we can suppose, after
suitable base extension, that the $2$-torsion is rational. Letting $B$
be the new base, the elliptic scheme can be viewed as a pull-back of the
Legendre scheme via a finite map $B\to B_0$. The height of the fiber is
then calculated by looking at the fibers over $B_0$.

{\bf (iii) Height of a section $\sigma: B \to \cE$. } 
 We follow \cite{cdmz}, where we gave a closed integral formula for the
 height of sections.
Given an elliptic scheme $\pi: \cE\to B$, and a simply connected open
 set $U\subset B$, we consider a basis $\rho_1,\rho_2$ of the period
 lattice: we can view them as sections of the line bundle
 $\Lie(\cE)$ over $U$; each point of $\cE$ over some point of
 $U$ admits a logarithm, defined up to integral multiples of
 $\rho_1,\rho_2$. This logarithm can be expressed as a linear
 combination with real coefficients of $\rho_1,\rho_2$. We then obtain
 real valued functions $\beta_1,\beta_2$ on $\pi^{-1}(U)$, which are
 locally well-defined, and globally  defined only up to addition of
 integer numbers. The differentials ${d}\beta_1,\,
 {d}\beta_2$ are well defined on the whole of
 $\pi^{-1}(U)$. Finally, the exterior product
\begin{equation*}
\omega:= {d}\beta_1 \wedge {d}\beta_2
\end{equation*}
 is well-defined on the whole of $\cE$. Note that its integral on each fiber equals $1$. 
 Now, for a rational section $\sigma:B\to \cE$ we define its height to be
 \begin{equation}\label{E.Height-Integral}
 \hat{h}(\sigma)=\int_B \sigma^* \omega.
 \end{equation}
 Clearly, it vanishes if and only if a non-trivial linear combination of
 $\beta_1,\beta_2$ is constant. By a theorem of Manin (see
 \cite{cz-poncelet} for a modern presentation and \cite{acz} for
 generalizations), this happens if and only if
$\sigma$ is torsion, or the scheme is isotrivial and the section is constant.

In \cite{cdmz}, \S 3 and \S 4, the following properties  have been established
\begin{enumerate}
\item The differential $2$-form $\omega$ is the only closed $2$-form on
      the total space $\cE$ such that $[2]^*\omega = 4\cdot \omega$,
      where $[2]:\cE\to\cE$ is the multiplication-by-2 morphism.
\item It is also the only $2$-form on $\cE$ whose integral on each fiber
      equals $1$, and such that $[2]^*\omega = 4\cdot \omega$.
\item The height $\hat{h}(\sigma)$ coincides with the normalized
      N\'eron-tate height as defined e.g. in chap. VI of Silverman's
      book \cite{silver}:
\begin{equation*}
\hat{h}(\sigma)=\lim_{n\to\infty} \frac{(n\sigma\cdot O)}{n^2},
\end{equation*}
where the numerator denotes the intersection product (on a smooth
      projective model of $\cE$) between the closure of the image of the
      section $n\sigma:B\to \cE$ and the zero section.
\item The height is also expressible as the limit
\begin{equation*}
\lim_{n\to\infty} \frac{\sharp \{b\in B\, |\, n\sigma(b)=0\}}{n^2}.
\end{equation*}
\item In the case of the Legendre scheme, the differential form $\omega$
      can be expressed as
\begin{equation*}
\omega =    dd^c (\Re(z\eta_\lambda(z)),
\end{equation*}
where $\lambda\in B$ is the coordinate on the base, $z$ is the
      coordinate in the Lie algebra $\C$ of $\cE_\lambda$, and the
      function $(\lambda,z)\mapsto \eta_\lambda(z)$ is defined as
      follows: (1) for each $\lambda\in B$, the function $z\mapsto
      \eta_\lambda(z)$ is $\R$-linear; (2) for each $\lambda$ with
      $|\lambda|<1, |\lambda-1|<1$, so that $\rho_1(\lambda),
      \rho_2(\lambda)$ are well defined by \eqref{E.Gauss},
      $\eta_\lambda(\rho_1(\lambda))$ and
      $\eta_\lambda(\rho_2(\lambda))$ are the semi-periods of
      $\cE_\lambda$; (3) the holomorphic functions $\lambda\mapsto
      \eta_\lambda(\rho_i(\lambda))$ are analytically continued along
      paths on $B$ (see again chap. VI of \cite{silver}).
\end{enumerate}

Whenever $\sigma: B\to \cE$ is a holomorphic section in general,
 the integral \eqref{E.Height-Integral} does not 
necessarily converge. 
In \cite{cz}, \S 3.2, we gave another integral formula for
 $\hat{h}(\sigma)$, this time in terms of the `modular logarithm',
 defined therein.

In order to extend the notion of
 normalized (N\'eron) height to transcendental sections,
 we follow the same pattern as in  \S \ref{rmkalg}:
Let us choose a rational function $\xi$ on $B$ such that $\xi$ is holomorphic on
$B$ and has a pole at every point of $\bar B \setminus B$, where
$\bar B$ is a smooth compactification of $B$.
Then, the modulus $|\xi|: B \to \R$ is a non-negative exhaustion function
 such that $dd^c \log |\xi| \equiv 0$ on $\{x \in B: \xi(x)\not=0\}$.
Setting $B(r)=\{x \in B: |\xi(x)|<r\}$ for $r>0$, we define the
{\em height (characteristic) function} by
\begin{equation*}
\hat{T}_\sigma(r) = \int_{0}^r \frac{{d}s}{s} \int_{B(s)} \sigma^*\omega.
\end{equation*}

In view of Theorem  \ref{rat}  and property (iv)
of the normalized height, whenever the section $\sigma$ is
transcendental, the integral \eqref{E.Height-Integral}
diverges. Vice-versa, whenever the height function   satisfies
$\hat{T}_\sigma(r)=O(\log r)$, then the integral \eqref{E.Height-Integral}
 converges, and again by
combining Theorem \ref{rat} with property (iv) one can deduce that
$\sigma$ is rational.

\bigskip
\begin{flushright}
Pietro Corvaja\\
Dipartimento di Scienze Matematiche, Informatiche e Fisiche\\
 Universit\`a di Udine, Via delle Scienze, 206, Udine, Italy\\
E-mail address: pietro.corvaja@uniud.it\\ \smallskip
Junjiro Noguchi\\
Graduate School of Mathematical Scineces, University of Tokyo\\
Komaba 3-8-1, Meguro-ku, Tokyo 153-8914, Japan\\
E-mail address: noguchi@ms.u-tokyo.ac.jp\\ \smallskip
Umberto Zannier\\
Scuola Normale Superiore, Piazza dei Cavalieri, 7\\
56126 Pisa, Italy\\
E-mail address: u.zannier@sns.it
\end{flushright}
\end{document}